\documentstyle[amstex,makeidx,12pt]{article}

\newtheorem{th}{Theorem}[section]
\newtheorem{lem}[th]{Lemma}
\newtheorem{prop}[th]{Proposition}
\newtheorem{cor}[th]{Corollary}
\newtheorem{defn}[th]{Definition}
\newenvironment{defn-new}{\begin{defn} \em}{\end{defn}}
\newtheorem{rem}[th]{Remark}
\newenvironment{rem-new}{\begin{rem} \em}{\end{rem}}
\newtheorem{ex}[th]{Example}
\newenvironment{ex-new}{\begin{ex} \em}{\end{ex}}
\newtheorem{exer}[th]{Exercise}
\newenvironment{exer-new}{\begin{exer} \em}{\end{exer}}
\newtheorem{agr}[th]{Agreement}
\newenvironment{agr-new}{\begin{agr} \em}{\end{agr}}
\newtheorem{pbm}[th]{Problem}
\newenvironment{pbm-new}{\begin{pbm} \em}{\end{pbm}}

\makeatletter \@addtoreset{equation}{section} \makeatother

\begin{document}

\begin{center}
{\Large {\bf Biharmonic Curves in $3$-dimensional Hyperbolic Heisenberg Group%
}}

\medskip

Selcen Y\"{u}ksel Perkta\c{s}, Erol K\i l\i \c{c}
\end{center}

\noindent {\bf Abstract.} In this paper we study the non-geodesic non-null
biharmonic curves in $3$-dimensional hyperbolic Heisenberg group. We prove
that all of the non-geodesic non-null biharmonic curves in $3$-dimensional
hyperbolic Heisenberg group are helices. Moreover, we obtain explicit
parametric equations for non-geodesic non-null biharmonic curves and
non-geodesic spacelike horizontal biharmonic curves in $3$-dimensional
hyperbolic Heisenberg group, respectively. We also show that there do not
exist non-geodesic timelike horizontal biharmonic curves in $3$-dimensional
hyperbolic Heisenberg group.



\section{Introduction\label{sect-intro}}

In 1964, Eells and Sampson \cite{Eells-Sampson} introduced the notion of
biharmonic maps as a natural generalization of the well-known harmonic maps.
Thus, while a map $\Psi $ from a compact Riemannian manifold $(M,g)$ to
another Riemannian manifold $(N,h)$ is harmonic if it is a critical point of
the energy functional $E(\Psi )=\frac{1}{2}\int_{M}|d\Psi |^{2}v_{g}$, the
biharmonic maps are the critical points of the bienergy functional $%
E_{2}(\Psi )=\frac{1}{2}\int_{M}|\tau (\Psi )|^{2}v_{g}.$

In a different setting, Chen \cite{Chen} defined biharmonic submanifolds $%
M\subset E^{n}$ of the Euclidean space as those with harmonic mean curvature
vector field, that is $\Delta H=0,$ where $\Delta $ is the rough Laplacian,
and stated that any biharmonic submanifold of the Euclidean space is
harmonic, that is minimal.

If the definition \ of biharmonic maps is applied to Riemannian immersions
into Euclidean space, the notion of Chen's biharmonic submanifold is
obtained, so the two definitions agree.

Harmonic maps are characterized by the vanishing of the tension field $\tau
(\Psi )=trace\nabla d\Psi $ where $\nabla $ is a connection induced from the
Levi-Civita connection $\nabla ^{M}$ of $M$ and $\nabla ^{\Psi }$ is the
pull-back connection. The first variation formula for the bienergy derived
in \cite{Jiang1, Jiang2} shows that the Euler-Lagrange equation for the
bienergy is
\[
\tau _{2}(\Psi )=-J(\tau (\Psi ))=-\Delta \tau (\Psi )-traceR^{N}(d\Psi
,\tau (\Psi ))d\Psi =0,
\]%
where $\Delta =-trace(\nabla ^{\Psi }\nabla ^{\Psi }-\nabla _{\nabla }^{\Psi
})$ is the rough Laplacian on the sections of $\Psi ^{-1}TN$ and $%
R^{N}(X,Y)=[\nabla _{X},\nabla _{Y}]-\nabla _{\lbrack X,Y]\text{ }}$ is the
curvature operator on $N$. From the expression of the bitension field $\tau
_{2}$, it is clear that a harmonic map is automatically a biharmonic map.
Non-harmonic biharmonic maps are called proper biharmonic maps.

Of course, the first and easiest examples can be found by looking at
differentiable curves in a Riemannian manifold. Obviously geodesics are
biharmonic. So, non-geodesic biharmonic curves are \ more interesting. Chen
and Ishikawa \cite{D13} showed non-existence of proper biharmonic curves in
Euclidean 3-space $E^{3}.$ Moreover they classified all proper biharmonic
curves in Minkowski 3-space $E_{1}^{3}$ (see also \cite{A12}). Caddeo,
Montaldo and Piu showed that on a surface with non-positive Gaussian
curvature, any biharmonic curve is a geodesic of the surface \cite{C14}. So
they gave a positive answer to generalized Chen's conjecture. Caddeo et al.
in \cite{A2} studied biharmonic curves in the unit 3-sphere. More precisely,
they showed that proper biharmonic curves in $S^{3}$ are circles of geodesic
curvature 1 or helices which are geodesics in the Clifford minimal torus.

On the other hand, there are several classification results on biharmonic
curves in arbitrary Riemannian manifolds. The biharmonic curves in the
Heisenberg group $H_{3}$ are investigated in \cite{D10} by Caddeo et al.
They showed that biharmonic curves in\ $H_{3}$ are helices, that is curves
with constant geodesic curvature $k_{1}$ and geodesic torsion $k_{2}.$ The
authors in \cite{Essin-Turhan-Arabian} studied non-geodesic horizontal
biharmonic curves in $3$-dimensional Heisenberg group. In \cite{A10} Fetcu
studied biharmonic curves in the generalized Heisenberg group and obtained
two families of proper biharmonic curves. Also, the explicit parametric
equations for the biharmonic curves on Berger spheres $S_{\varepsilon }^{3}$
are obtained by Balmu\c{s} in \cite{C6}.

In contact geometry, there is a well known analog of real space form, namely
a Sasakian space form. In particular, a simply connected three-dimensional
Sasakian space form of constant holomorphic sectional curvature $1$ is
isometric to $S^{3}.$ So in this context J. Inoguchi classified in \cite{A11}
the proper biharmonic Legendre curves and Hopf cylinders in a $3$%
-dimensional Sasakian space form and in \cite{B11} the explicit parametric
equations were obtained. In \cite{Cho-Inoguchi-Lee}, the authors showed that
every non-geodesic biharmonic curve in a $3$-dimensional Sasakian space form
of constant holomorphic sectional curvature is a helix. T. Sasahara \cite%
{C53}, analyzed the proper biharmonic Legendre surfaces in Sasakian space
forms and in the case when the ambient space is the unit $5$-dimensional
sphere $S^{5}$ he obtained their explicit representations. A full
classification of proper biharmonic Legendre curves, explicit examples and a
method to construct proper biharmonic anti-invariant submanifolds in any
dimensional Sasakian space form were given in \cite{Fetcu-Oniciuc-2009}.
Furthermore, D. Fetcu \cite{Fetcu-2010} studied proper biharmonic
non-Legendre curves in a Sasakian space form.

Geometry of almost paracontact manifolds can be considered as a natural
extension of the almost paraHermitian geometry to the odd dimensional case
while the almost contact manifolds are a natural extension of the almost
Hermitian manifolds. A paracontact structure on a real $(2n+1)$-dimensional
manifold $M$ is a $(1,1)$ tensor field $\varphi ,$ a vector field $\xi $, a
codimension one distribution $D$ (horizontal bundle), a paracomplex
structure $I|_{\varphi }$ on $D$, that is, $I^{2}=id$ and the $\pm $
eigendistributions $D^{\pm }$ have equal dimension. Locally, the horizontal
bundle $D$ is given by the kernel of a $1$-form $\eta ,$ that is $D=\ker
\eta $. A paracontact structure is called a paracontact Hermitian structure
if $\eta $ is a para Hermitian contact form in the sense that there exist a
non-degenerate semi-Riemannian metric $g$, which is defined on $D$, and
compatible with $\eta $ and $I$, $d\eta (X,Y)=2g(IX,Y)$, \ $%
g(IX,IY)=-g(X,Y), $ for all $X,$ $Y\in D.$ The signature of $g$ on $D$ is
necessarily of (signature) type $(n,n)$. If the paracomplex structure $I$ on
$D$ is integrable, that is $[D^{\pm },D^{\pm }]\in D^{\pm }$, then the
paracontact structure is said to be integrable. A paracontact manifold with
an integrable paracontact structure is called a para CR-manifold. A
paracontact manifold is said to be paraSasakian if $N(X,Y)=2d\eta (X,Y)\xi $%
, where $N$ is the Nijenhuis tensor of $I$ given by $%
N(X,Y)=[IX,IY]+[X,Y]-I[IX,Y]-I[X,IY] $, $X,$ $Y\in D$ \cite%
{Zamkovoy-Global-2009}, \cite{Ivanov-Vassilev-Zamkovoy-Dedicate-2010}.

The basic example of a paracontact manifold is the hyperbolic Heisenberg
group. Let $G(P)=R^{2n}\times R$ \ be a group with the group law given by%
\[
(p^{\prime },t^{\prime })\circ \left( p,t\right) =\left( p^{\prime
}+p,t^{\prime }+t-\sum_{k=1}^{n}\left( u_{k}^{\prime }v_{k}-v_{k}^{\prime
}u_{k}\right) \right)
\]%
where $p^{\prime }=\left( u_{1}^{\prime },v_{1}^{\prime },...,u_{n}^{\prime
},v_{n}^{\prime }\right) ,$ $p=\left( u_{1},v_{1},...,u_{n},v_{n}\right) \in
R^{2n}$ and $t^{\prime },t\in R.$ A basis of left-invariant vector fields is
given by
\[
U_{k}=\frac{\partial }{\partial u_{k}}-2v_{k}\frac{\partial }{\partial t}%
,\quad V_{k}=\frac{\partial }{\partial v_{k}}-2u_{k}\frac{\partial }{%
\partial t},\quad \xi =2\frac{\partial }{\partial t}.
\]%
Define $\widetilde{\Theta }=-\frac{1}{2}dt-\sum_{k=1}^{n}\left(
u_{k}dv_{k}-v_{k}du_{k}\right) $ with corresponding horizontal distribution $%
D$ given by the span of the left invariant horizontal vector fields $\left\{
U_{1},...,U_{n},V_{1},...,V_{n}\right\} .$ An endomorphism on $D$ defined by
$IU_{k}=V_{k},$ $IV_{k}=U_{k}$ is a paracomplex structure on $D.$ The form $%
\widetilde{\Theta }$ and the paracomplex structure $I$ define a paracontact
manifold which is called the hyperbolic Heisenberg group and denoted by $%
\left( G(P),\eta \right) .$ Note that $\left\{
U_{1},...,U_{n},V_{1},...,V_{n},\xi \right\} $ is an orthonormal basis of
the tangent space, $g(U_{j},U_{j})=-g(V_{j},V_{j})=1,$ $1\leq j\leq n$ \cite%
{Zamkovoy-Global-2009}, \cite{Ivanov-Vassilev-Zamkovoy-Dedicate-2010}. The
authors in \cite{Ivanov-Vassilev-Zamkovoy-Dedicate-2010} also proved that an
integrable paracontact Hermitian manifold $(M,\eta ,I,g)$ of dimension $3$
is locally isomorphic to the $3$-dimensional hyperbolic Heisenberg group
exactly when the canonical connection has vanishing horizontal curvature and
zero torsion. So, this motivated us to initiate study of the biharmonic
curves in paracontact manifolds by studying biharmonic curves in $3$%
-dimensional hyperbolic Heisenberg group.

In this paper we study the non-null biharmonic curves in $3$-dimensional
hyperbolic Heisenberg group (for short, ${\cal HH}_{3}$). Section 1 is
devoted to the some basic definitions. We also define and characterize a
cross product in $3$-dimensional hyperbolic Heisenberg group. In section 2
we investigate the necessary and sufficient conditions for a non-null curve
in $3$-dimensional hyperbolic Heisenberg group to be non-geodesic
biharmonic. In section 3 we prove that a non-geodesic non-null curve
parametrized by arclenght in $3$-dimensional hyperbolic Heisenberg group
with the vanishing third component of the binormal vector field cannot be
biharmonic. In section 4, we study the non-geodesic non-null biharmonic
helices in $3$-dimensional hyperbolic Heisenberg group. Moreover, we obtain
explicit parametric equations for non-geodesic non-null biharmonic curves in
$3$-dimensional hyperbolic Heisenberg group. In the last section, we give
explicit examples of non-geodesic spacelike horizontal biharmonic curves and
prove that there do not exist non-geodesic timelike horizontal biharmonic
curves in $3$-dimensional hyperbolic Heisenberg group.

\section{Preliminaries\label{pre}}

\subsection{Biharmonic Maps}

Let $(M,g)$ and $(N,h)$ be Riemannian manifolds and $\Psi :(M,g)\rightarrow
(N,h)$ \ be a smooth map. The tension field of $\Psi $ is given by $\tau
(\Psi )=trace\nabla d\Psi $, where $\nabla d\Psi $ is the second fundamental
form of $\Psi $ defined by $\nabla d\Psi (X,Y)=\nabla _{X}^{\Psi }d\Psi
(Y)-d\Psi (\nabla _{X}^{M}Y)$, $X,\,Y\in \Gamma (TM)$. For any compact
domain $\Omega \subseteq M$, the bienergy is defined by
\[
E_{2}(\Psi )=\frac{1}{2}\int_{\Omega }|\tau (\Psi )|^{2}v_{g}.
\]%
Then a smooth map $\Psi $ is called biharmonic map if it is a critical point
of the bienergy functional for any compact domain $\Omega \subseteq M.$ We
have for the bienergy the following first variation formula:
\[
\frac{d}{dt}E_{2}(\Psi _{t};\Omega )|_{t=0}=\text{ }\int_{\Omega }<\tau
_{2}(\Psi ),w>v_{g}
\]%
where $v_{g}$ is the volume element, $w$ is the variational vector field
associated to the variation $\{\Psi _{t}\}$ of $\Psi $ and
\[
\tau _{2}(\Psi )=-J(\tau _{2}(\Psi ))=-\Delta ^{\Psi }\tau (\Psi
)-traceR^{N}(d\Psi ,\tau (\Psi ))d\Psi .
\]%
$\tau _{2}(\Psi )$ is called bitension field of $\Psi $. Here $\Delta ^{\Psi
}$ is the rough Laplacian on the sections of the pull-back bundle $\Psi
^{-1}TN$ which is defined by
\[
\Delta ^{\Psi }V=-\sum_{i=1}^{m}\{\nabla _{e_{i}}^{\Psi }\nabla
_{e_{i}}^{\Psi }V-\nabla _{\nabla _{e_{i}}^{M}e_{i}}^{\Psi
}V\},\,\,\,\,\,V\in \Gamma (\Psi ^{-1}TN),
\]%
where $\nabla ^{\Psi }$ is the pull-back connection on the pull-back bundle $%
\Psi ^{-1}TN$ and $\{e_{i}\}_{i=1}^{m}$ is an orthonormal frame on $M.$ When
the target manifold is semi-Riemannian manifold, the bienergy and bitension
field can be defined in the same way.

Let $M$ be a semi-Riemannian manifold and $\gamma :I\rightarrow {M}$ be a
non-null curve parametrized by arclenght. By using the definition of the
tension field we have
\[
\tau (\gamma )=\nabla _{\frac{\partial }{\partial s}}^{\gamma }d\gamma (%
\frac{\partial }{\partial s})=\nabla _{T}T,
\]%
where $T=\gamma ^{\prime }$. In this case biharmonic equation for the curve $%
\gamma $ reduces to
\[
\tau _{2}(\gamma )=\nabla _{T}^{3}T-R(T,\nabla _{T}T)T=0.
\]

\subsection{3-dimensional Hyperbolic Heisenberg Group}

Consider $R^{3}$ with the group law given by
\begin{equation}
\widetilde{X}X=(\widetilde{x}+x,\widetilde{y}+y,\widetilde{z}+z-\widetilde{x}%
y+\widetilde{y}x),  \label{hh1}
\end{equation}%
where $X=(x,y,z),$ $\widetilde{X}=(\widetilde{x},\widetilde{y},\widetilde{z}%
) $.

Let ${\cal HH}_{3}=(R^{3},g)$ be $3$-dimensional hyperbolic Heisenberg group
endowed with the semi-Riemannian metric $g$ which is defined by
\begin{equation}
g=(dx)^{2}+(dy)^{2}-\frac{1}{4}(dz+2ydx-2xdy)^{2}.  \label{hh2}
\end{equation}%
Note that the metric $g$ is left invariant.

We can define an orthonormal basis for the tangent space of ${\cal HH}_{3}$
by
\begin{equation}
e_{1}=\frac{\partial }{\partial x}-2y\frac{\partial }{\partial z},\quad
e_{2}=\frac{\partial }{\partial y}+2x\frac{\partial }{\partial z},\quad
e_{3}=2\frac{\partial }{\partial z},  \label{hh3}
\end{equation}%
which is dual to the coframe%
\[
\theta ^{1}=dx,\quad \theta ^{2}=dy,\quad \theta ^{3}=\frac{1}{2}dz+ydx-xdy.
\]

\begin{prop}
For the covariant derivatives of the Levi-Civita connection of the
left-invariant metric $g$ defined above, we have
\begin{equation}
\left\{
\begin{array}{c}
\nabla _{e_{1}}e_{1}=0,\quad \nabla _{e_{1}}e_{2}=e_{3},\quad \nabla
_{e_{1}}e_{3}=-e_{2}, \\
\nabla _{e_{2}}e_{1}=-e_{3},\quad \nabla _{e_{2}}e_{2}=0,\quad \nabla
_{e_{2}}e_{3}=-e_{1}, \\
\nabla _{e_{3}}e_{1}=-e_{2},\quad \nabla _{e_{3}}e_{2}=-e_{1},\quad \nabla
_{e_{3}}e_{3}=0,%
\end{array}%
\right.  \label{hh4}
\end{equation}%
where $\{e_{1},e_{2},e_{3}\}$ is the orthonormal basis for the tangent space
given (\ref{hh3})
\end{prop}

Also, we have the following bracket relations%
\begin{equation}
\left[ e_{1},e_{2}\right] =2e_{3},\quad \left[ e_{1},e_{3}\right] =\left[
e_{2},e_{3}\right] =0.  \label{hh5}
\end{equation}%
The curvature tensor field of $\nabla $ is given by
\[
R(X,Y)Z=\nabla _{X}\nabla _{Y}Z-\nabla _{Y}\nabla _{X}Z-\nabla _{\left[ X,Y%
\right] }Z,
\]%
while the Riemannian-Christoffel tensor field is
\[
R(X,Y,Z,W)=g(R(X,Y)Z,W),
\]%
where $X,$ $Y,$ $Z,$ $W\in \Gamma \left( T{\cal H}_{3}\right) $. If we put
\[
R_{abc}=R(e_{a},e_{b})e_{c},
\]%
where the indices $a,$ $b,$ $c$ take the values $1,$ $2,$ $3$. Then the
non-zero components of the curvature tensor field are
\begin{equation}
\left\{
\begin{array}{c}
R_{121}=3e_{2},\quad R_{122}=3e_{1},\quad R_{131}=-e_{3},\quad \\
R_{133}=-e_{1},\quad R_{232}=e_{3},\quad R_{233}=-e_{2}.%
\end{array}%
\right.  \label{hh6}
\end{equation}

Now we shall define a cross product on $3$-dimensional hyperbolic Heisenberg
group for later use

\begin{defn-new}
We define a cross product $\wedge $ on ${\cal HH}_{3}$ by
\[
X\wedge Y=-\left( a_{2}b_{3}-a_{3}b_{2}\right) e_{1}-\left(
a_{1}b_{3}-a_{3}b_{1}\right) e_{2}+\left( a_{1}b_{2}-a_{2}b_{1}\right)
e_{3},
\]%
where $\{e_{1},e_{2},e_{3}\}$ is an orthonormal basis of ${\cal HH}_{3}$
given by (\ref{hh3}) and $X=a_{1}e_{1}+a_{2}e_{2}+a_{3}e_{3},$ $%
Y=b_{1}e_{1}+b_{2}e_{2}+b_{3}e_{3}\in \Gamma \left( T({\cal HH}_{3})\right) $%
.
\end{defn-new}

\begin{th}
The cross product $\wedge $ on ${\cal HH}_{3}$ has the following properties:

\begin{description}
\item[(i)] The cross product is bilinear and anti-symmetric ($X\wedge
Y=-Y\wedge X$).

\item[(ii)] $X\wedge Y$ is perpendicular both of $X$ and $Y$.

\item[(iii)] $e_{1}\wedge e_{2}=e_{3},$\quad $e_{2}\wedge e_{3}=-e_{1},$%
\quad $e_{3}\wedge e_{1}=e_{2}.$

\item[(iv)] $\left( X\wedge Y\right) \wedge Z=g(X,Z)Y-g(Y,Z)X.$

\item[(v)] Define a mixed product by%
\[
\left( X,Y,Z\right) =g(X\wedge Y,Z),
\]%
then we have
\[
\left( X,Y,Z\right) =-\det (X,Y,Z)
\]%
and
\[
\left( X,Y,Z\right) =(Y,Z,X)=(Z,X,Y).
\]

\item[(vi)] $\left( X\wedge Y\right) \wedge Z+\left( Y\wedge Z\right) \wedge
X+\left( Z\wedge X\right) \wedge Y=0,$
\end{description}
\end{th}

for all $X,$ $Y,$ $Z\in \Gamma \left( T({\cal HH}_{3})\right) .$

\section{Biharmonic curves in $3$-dimensional hyperbolic\newline
Heisenberg group}

An arbitrary curve $\gamma :I\rightarrow {\cal HH}_{3},$ $\gamma =\gamma
(s), $ in $3$-dimensional hyperbolic Heisenberg group ${\cal HH}_{3}$ is
called spacelike, timelike or null (lightlike), if all of its velocity
vectors $\gamma ^{\prime }(s)$ are respectively spacelike, timelike or null
(lightlike). If $\gamma (s)$ is a spacelike or timelike curve, we can
reparametrize it such that $g(\gamma ^{\prime }(s),\gamma ^{\prime
}(s))=\varepsilon $ where $\varepsilon =1$ if $\gamma $ is spacelike and $%
\varepsilon =-1$ if $\gamma $ is timelike, respectively. In this case $%
\gamma (s)$ is said to be unit speed or arclenght parametrization.

Let $\gamma :I\rightarrow {\cal HH}_{3}$ be a non-null curve parametrized by
arclenght and $\{T,N,B\}$ be the orthonormal moving Frenet frame along the
curve $\gamma $ in ${\cal HH}_{3}$ such that $T=\gamma ^{^{\prime }}$ is the
unit vector field tangent to $\gamma ,$ $N$ is the unit vector field in the
direction $\nabla _{T}T$ normal to $\gamma $ and $B=T\wedge N$. The mutually
orthogonal unit vector fields $T$, $N$ and $B$ are called the tangent, the
principal normal and the binormal vector fields, respectively. Then we have
the following Frenet equations%
\begin{eqnarray}
\nabla _{T}T &=&k_{1}\varepsilon _{2}N,  \nonumber \\
\nabla _{T}N &=&-k_{1}\varepsilon _{1}T+k_{2}\varepsilon _{3}B,
\label{bhc-1} \\
\nabla _{T}B &=&-k_{2}\varepsilon _{2}N,  \nonumber
\end{eqnarray}%
where $\varepsilon _{1}=g(T,T),$ $\varepsilon _{2}=g(N,N)$ and $\varepsilon
_{3}=g(B,B).$ Here $k_{1}=\left\vert \tau (\gamma )\right\vert =\left\vert
\nabla _{T}T\right\vert $ is the geodesic curvature of $\gamma $ and $k_{2}$
is its geodesic torsion.

From (\ref{bhc-1}) we have
\begin{eqnarray}
\nabla _{T}^{3}T &=&\left( -3k_{1}k_{1}^{\prime }\varepsilon _{1}\varepsilon
_{2}\right) T+\left( k_{1}^{\prime \prime }\varepsilon
_{2}-k_{1}^{3}\varepsilon _{1}-k_{1}k_{2}^{2}\varepsilon _{3}\right) N
\nonumber \\
&&+\left( 2k_{1}^{\prime }k_{2}\varepsilon _{2}\varepsilon
_{3}+k_{1}k_{2}^{\prime }\varepsilon _{2}\varepsilon _{3}\right) B.
\label{bhc-2}
\end{eqnarray}%
Using (\ref{hh6}) one obtains
\begin{subequations}
\begin{equation}
R(T,\nabla _{T}T)T=k_{1}\varepsilon _{2}\left[ \left( -\varepsilon
_{2}\varepsilon _{3}-4\varepsilon _{2}B_{3}^{2}\right) N+\left( 4\varepsilon
_{3}N_{3}B_{3}\right) B\right] ,  \label{bhc-3}
\end{equation}%
where $T=T_{1}e_{1}+T_{2}e_{2}+T_{3}e_{3},$ $%
N=N_{1}e_{1}+N_{2}e_{2}+N_{3}e_{3}$ and $B=T\wedge
N=B_{1}e_{1}+B_{2}e_{2}+B_{3}e_{3}$. Hence we get
\end{subequations}
\begin{eqnarray}
\tau _{2}(\gamma ) &=&\left( -3k_{1}k_{1}^{\prime }\varepsilon
_{1}\varepsilon _{2}\right) T+\left( k_{1}^{\prime \prime }\varepsilon
_{2}-k_{1}^{3}\varepsilon _{1}-k_{1}k_{2}^{2}\varepsilon
_{3}+k_{1}\varepsilon _{3}+4k_{1}B_{3}^{2}\right) N  \nonumber \\
&&+\left( 2k_{1}^{\prime }k_{2}\varepsilon _{2}\varepsilon
_{3}+k_{1}k_{2}^{\prime }\varepsilon _{2}\varepsilon _{3}-4k_{1}\varepsilon
_{2}\varepsilon _{3}N_{3}B_{3}\right) B.  \label{bhc-4}
\end{eqnarray}

\begin{th}
Let $\gamma :I\rightarrow {\cal HH}_{3}$ be a non-null curve parametrized by
arclenght. Then $\gamma $ is a non-geodesic biharmonic curve if and only if
\begin{equation}
\left\{
\begin{array}{c}
k_{1}=\text{constant}\neq 0, \\
k_{1}^{2}\varepsilon _{1}\varepsilon _{3}+k_{2}^{2}=1+4\varepsilon
_{3}B_{3}^{2}, \\
k_{2}^{\prime }=N_{3}B_{3}.%
\end{array}%
\right.  \label{bhc-5}
\end{equation}
\end{th}

\noindent {\em Proof.} From (\ref{bhc-4}) it follows that $\gamma $ is
biharmonic if and only if%
\[
\left\{
\begin{array}{c}
k_{1}k_{1}^{\prime }=0, \\
k_{1}^{\prime \prime }\varepsilon _{2}-k_{1}^{3}\varepsilon
_{1}-k_{1}k_{2}^{2}\varepsilon _{3}+k_{1}\varepsilon _{3}+4k_{1}B_{3}^{2}=0,
\\
2k_{1}^{\prime }k_{2}+k_{1}k_{2}^{\prime }-4k_{1}N_{3}B_{3}=0.%
\end{array}%
\right.
\]%
If we look for non-geodesic solution of the above system we complete the
proof.

\begin{cor}
If $k_{1}=$constant$\neq 0$ and $k_{2}=0$ for a non-null curve $\gamma
:I\rightarrow {\cal HH}_{3}$ then $\gamma $ is a non-geodesic biharmonic
curve if and only if $k_{1}^{2}=\varepsilon _{1}\left( \varepsilon
_{3}+4B_{3}^{2}\right) $ and $N_{3}B_{3}=0.$
\end{cor}

\begin{prop}
Let $\gamma :I\rightarrow {\cal HH}_{3}$ be a non-geodesic non-null curve
parametrized by arclenght. If $k_{1}$ is constant and $N_{3}B_{3}\neq 0$,
then $\gamma $ is not biharmonic.
\end{prop}

\noindent {\em Proof.} By using (\ref{hh4}) and (\ref{bhc-1}) we have
\begin{eqnarray}
\nabla _{T}T &=&\left( T_{1}^{\,\prime }-2T_{2}T_{3}\right) e_{1}+\left(
T_{2}^{\,\prime }-2T_{1}T_{3}\right) e_{2}+T_{3}^{\,\prime }e_{3}
\label{bhc-7} \\
&=&k_{1}\varepsilon _{2}N,  \nonumber
\end{eqnarray}%
which implies that
\[
T_{3}^{\,\prime }=k_{1}\varepsilon _{2}N_{3}.
\]%
If we put $T_{3}(s)=k_{1}F(s)$ and $f(s)=F^{\prime }(s)$ we get $%
f(s)=\varepsilon _{2}N_{3}(s).$ Then we can write
\[
T=\sqrt{\varepsilon _{1}+k_{1}^{2}F^{2}}\cosh \beta \,e_{1}+\sqrt{%
\varepsilon _{1}+k_{1}^{2}F^{2}}\sinh \beta \,e_{2}+k_{1}Fe_{3}.
\]%
From (\ref{bhc-7}) we calculate
\begin{eqnarray}
\nabla _{T}T &=&k_{1}\varepsilon _{2}N=\left( \frac{k_{1}^{2}Ff}{\sqrt{%
\varepsilon _{1}+k_{1}^{2}F^{2}}}\cosh \beta +\sqrt{\varepsilon
_{1}+k_{1}^{2}F^{2}}(\beta ^{\prime }-k_{1}F)\sinh \beta \right) e_{1}
\nonumber \\
&&\quad \qquad \quad +\left( \frac{k_{1}^{2}Ff}{\sqrt{\varepsilon
_{1}+k_{1}^{2}F^{2}}}\sinh \beta +\sqrt{\varepsilon _{1}+k_{1}^{2}F^{2}}%
(\beta ^{\prime }-k_{1}F)\cosh \beta \right) e_{2}  \nonumber \\
&&\qquad \qquad +\left( k_{1}\,f\right) e_{3}.  \label{bhc-8}
\end{eqnarray}%
By taking into account the definition of the geodesic curvature $k_{1}$ and
the last equation one can see that
\begin{equation}
\beta ^{\prime }-k_{1}F=\pm k_{1}\frac{\sqrt{-\varepsilon _{1}\varepsilon
_{2}-\varepsilon _{1}f^{\,2}-\varepsilon _{2}k_{1}^{2}F^{2}}}{\varepsilon
_{1}+k_{1}^{2}F^{2}}.  \label{bhc-9}
\end{equation}%
If we write (\ref{bhc-9}) in (\ref{bhc-8}) we get%
\begin{eqnarray*}
\varepsilon _{2}N &=&\left( \pm \frac{\sqrt{-\varepsilon _{1}\varepsilon
_{2}-\varepsilon _{1}f^{\,2}-\varepsilon _{2}k_{1}^{2}F^{2}}}{\sqrt{%
\varepsilon _{1}+k_{1}^{2}F^{2}}}\sinh \beta +\frac{k_{1}Ff}{\sqrt{%
\varepsilon _{1}+k_{1}^{2}F^{2}}}\cosh \beta \right) e_{1} \\
&&+\left( \pm \frac{\sqrt{-\varepsilon _{1}\varepsilon _{2}-\varepsilon
_{1}f^{\,2}-\varepsilon _{2}k_{1}^{2}F^{2}}}{\sqrt{\varepsilon
_{1}+k_{1}^{2}F^{2}}}\cosh \beta +\frac{k_{1}Ff}{\sqrt{\varepsilon
_{1}+k_{1}^{2}F^{2}}}\sinh \beta \right) e_{2} \\
&&+f\,e_{3.}
\end{eqnarray*}%
Since $B=T\wedge N$, from the definition of the cross product in ${\cal H}%
_{3}$ we have
\begin{equation}
B_{3}=\pm \,\varepsilon _{2}\sqrt{-\varepsilon _{1}\varepsilon
_{2}-\varepsilon _{1}f^{\,2}-\varepsilon _{2}k_{1}^{2}F^{2}}.  \label{bhc-10}
\end{equation}%
On the other hand from the Frenet equations we obtain
\[
g(\nabla _{T}N,e_{3})=k_{1}\varepsilon _{1}T_{3}-k_{2}\varepsilon _{3}B_{3}.
\]%
Using (\ref{hh4}) since $N=N_{1}e_{1}+N_{2}e_{2}+N_{3}e_{3}$ we also have
\[
g(\nabla _{T}N,e_{3})=-N_{3}^{\,\prime }-B_{3},
\]%
which implies that
\begin{equation}
-N_{3}^{\,\prime }-B_{3}=k_{1}\varepsilon _{1}T_{3}-k_{2}\varepsilon
_{3}B_{3}.  \label{bhc-11}
\end{equation}%
By writing $N_{3}=\varepsilon _{2}f,$ $T_{3}=k_{1}F$ and (\ref{bhc-10}) in (%
\ref{bhc-11}) we get%
\begin{equation}
k_{2}=\pm \,\frac{\left( f^{\,\prime }+k_{1}^{2}\varepsilon _{1}\varepsilon
_{2}F\right) \varepsilon _{3}}{\sqrt{-\varepsilon _{1}\varepsilon
_{2}-\varepsilon _{1}f^{\,2}-\varepsilon _{2}k_{1}^{2}F^{2}}}+\varepsilon
_{3}=\pm \,\varepsilon _{1}\varepsilon _{3}\frac{B_{3}^{\,\prime }}{N_{3}}%
+\varepsilon _{3}.  \label{bhc-12}
\end{equation}%
Now assume that $\gamma $ is biharmonic. Then from the third equation in (%
\ref{bhc-5}) we write $k_{2}^{\prime }=N_{3}B_{3}\neq 0$ which gives
\[
N_{3}=\frac{k_{2}^{\prime }}{B_{3}}.
\]%
By writing the last equation in (\ref{bhc-12}) and then by integrating we
obtain
\begin{equation}
k_{2}^{\,2}=\pm \,\varepsilon _{1}\varepsilon
_{3}B_{3}^{2}+2k_{2}\varepsilon _{3}+c,  \label{bhc-13}
\end{equation}%
where $c$ is a constant. Also, from the second equation in (\ref{bhc-5}) we
have%
\begin{equation}
\varepsilon _{1}\varepsilon _{3}B_{3}^{2}=k_{1}^{2}\frac{\varepsilon _{3}}{4}%
+k_{2}^{2}\frac{\varepsilon _{1}}{4}-\frac{\varepsilon _{1}}{4}.
\label{bhc-14}
\end{equation}%
By comparing (\ref{bhc-13}) and (\ref{bhc-14}) we get
\[
k_{2}^{\,2}\left( 4\mp \varepsilon _{1}\right) -8k_{2}\varepsilon _{3}=C,
\]%
where $C=\mp $\thinspace $\varepsilon _{1}\pm k_{1}^{2}\,\varepsilon _{3}+4c$
is a constant, which implies that $k_{2}$ is also a constant. Hence we
obtain a contradiction with the assumption $k_{2}^{\prime }\neq 0$. This
completes the proof.

\begin{th}
Let $\gamma :I\rightarrow {\cal HH}_{3}$ be a non-geodesic non-null curve
parametrized by arclenght. Then $\gamma $ is biharmonic if and only if%
\begin{equation}
\left\{
\begin{array}{c}
k_{1}=\text{constant}\neq 0, \\
k_{2}=\text{constant,} \\
N_{3}B_{3}=0, \\
k_{1}^{2}\varepsilon _{1}\varepsilon _{3}+k_{2}^{2}=1+4\varepsilon
_{3}B_{3}^{2}.%
\end{array}%
\right.  \label{bhc-15}
\end{equation}%
\noindent
\end{th}

\section{Biharmonic helices in $3$-dimensional hyperbolic\newline
Heisenberg group}

A non-null curve in a semi-Riemannian manifold having constant both geodesic
curvature and geodesic torsion is called helix. Now we shall investigate the
biharmonicity conditions of a helix in $3$-dimensional hyperbolic Heisenberg
group. For any helix in ${\cal HH}_{3}$, the system (\ref{bhc-5}) reduces to

\begin{equation}
\left\{
\begin{array}{c}
k_{1}^{2}\varepsilon _{1}\varepsilon _{3}+k_{2}^{2}=1+4\varepsilon
_{3}B_{3}^{2}, \\
N_{3}B_{3}=0,%
\end{array}%
\right.  \label{bhc-6}
\end{equation}%
which implies that $B_{3}$ must be a constant.

\begin{prop}
\label{prop-bhh1}Let $\gamma :I\rightarrow {\cal HH}_{3}$ be a non-geodesic
non-null curve parametrized by arclenght with $B_{3}=0.$ Then we have $%
\varepsilon _{1}=-\varepsilon _{2}$ and $B$ is a timelike vector field,
where $\varepsilon _{1}=g(T,T)$ and $\varepsilon _{2}=g(N,N)$.
\end{prop}

\noindent {\em Proof.} Assume that $\gamma :I\rightarrow {\cal HH}_{3}$ is a
non-geodesic non-null curve parametrized by arclenght and $\gamma ^{\,\prime
}(s)=T(s)$. If $\gamma $ is a spacelike curve then we can write
\begin{equation}
T=\cosh \alpha _{1}\cosh \beta _{1}e_{1}+\cosh \alpha _{1}\sinh \beta
_{1}e_{2}+\sinh \alpha _{1}e_{3}.  \label{bhh1}
\end{equation}%
where $\alpha _{1}=\alpha _{1}(s)$ and $\beta _{1}=\beta _{1}(s)$. From (\ref%
{hh4}) the covariant derivative of the unit tangent vector field $T$, of $%
\gamma $, is
\begin{eqnarray*}
\nabla _{T}T &=&\left( \alpha _{1}^{\prime }\sinh \alpha _{1}\cosh \beta
_{1}+\cosh \alpha _{1}\sinh \beta _{1}\left( \beta _{1}^{\prime }-2\sinh
\alpha _{1}\right) \right) e_{1} \\
&&+\left( \alpha _{1}^{\prime }\sinh \alpha _{1}\sinh \beta _{1}+\cosh
\alpha _{1}\cosh \beta _{1}\left( \beta _{1}^{\prime }-2\sinh \alpha
_{1}\right) \right) e_{2} \\
&&+\left( \alpha _{1}^{\prime }\cosh \alpha _{1}\right) e_{3} \\
&=&k_{1}\varepsilon _{2}N.
\end{eqnarray*}%
By using the definition of cross product in ${\cal HH}_{3}$ we also obtain%
\[
B_{3}=\frac{\cosh ^{2}\alpha _{1}\left( \beta _{1}^{\prime }-2\sinh \alpha
_{1}\right) \varepsilon _{2}}{k_{1}}.
\]%
Now let $B_{3}=0$. From the last equation above, since $\cosh \alpha
_{1}\neq 0$ then $\beta _{1}^{\prime }-2\sinh \alpha _{1}=0$. Thus we have%
\begin{equation}
\nabla _{T}T=\alpha _{1}^{\prime }\left( \sinh \alpha _{1}\cosh \beta
_{1}e_{1}+\sinh \alpha _{1}\sinh \beta _{1}e_{2}+\cosh \alpha
_{1}e_{3}\right) .  \label{bhh2}
\end{equation}%
We can assume that $\alpha _{1}^{\prime }\neq 0$ ( when $\alpha _{1}^{\prime
}=0$ then we have $\nabla _{T}T=0$, which implies that $\gamma $ is a
geodesic). Hence we get%
\begin{eqnarray}
k_{1}^{\,2}\varepsilon _{2} &=&g(\nabla _{T}T,\nabla _{T}T)  \nonumber \\
&=&\left( \alpha _{1}^{\prime }\right) ^{2}\left( \sinh ^{2}\alpha
_{1}-\cosh ^{2}\alpha _{1}\right)  \nonumber \\
&=&-\left( \alpha _{1}^{\prime }\right) ^{2}.  \label{bhh2a}
\end{eqnarray}%
If $N$ is spacelike then $k_{1}=0$ which is a contradiction.

By a similar way, for a timelike curve $\gamma $, its tangent vector field
can be expressed by%
\begin{equation}
T=\sinh \alpha _{2}\cosh \beta _{2}e_{1}+\sinh \alpha _{2}\sinh \beta
_{2}e_{2}+\cosh \alpha _{2}e_{3}.  \label{bhh3}
\end{equation}%
where $\alpha _{2}=\alpha _{2}(s)$ and $\beta _{2}=\beta _{2}(s)$. From (\ref%
{hh4}) we get
\begin{eqnarray*}
\nabla _{T}T &=&\left( \alpha _{2}^{\prime }\cosh \alpha _{2}\cosh \beta
_{2}+\sinh \alpha _{2}\sinh \beta _{2}\left( \beta _{2}^{\prime }-2\cosh
\alpha _{2}\right) \right) e_{1} \\
&&+\left( \alpha _{2}^{\prime }\cosh \alpha _{2}\sinh \beta _{2}+\sinh
\alpha _{2}\cosh \beta _{2}\left( \beta _{2}^{\prime }-2\cosh \alpha
_{2}\right) \right) e_{2} \\
&&+\left( \alpha _{2}^{\prime }\sinh \alpha _{2}\right) e_{3} \\
&=&k_{1}\varepsilon _{2}N.
\end{eqnarray*}%
Next, we have
\begin{eqnarray*}
B_{3} &=&T_{1}N_{2}-T_{2}N_{1} \\
&=&\frac{\sinh ^{2}\alpha _{2}\left( \beta _{2}^{\prime }-2\cosh \alpha
_{2}\right) }{k_{1}}\varepsilon _{2}.
\end{eqnarray*}%
Now assume that $B_{3}=0$. If $\sinh \alpha _{2}=0$ then $T=e_{3}$, that is,
$\gamma $ is a geodesic. So one must have
\[
\beta _{2}^{\prime }-2\cosh \alpha _{2}=0.
\]%
Thus we get%
\begin{equation}
\nabla _{T}T=\alpha _{2}^{\prime }\left( \cosh \alpha _{2}\cosh \beta
_{2}e_{1}+\cosh \alpha _{2}\sinh \beta _{2}e_{2}+\sinh \alpha
_{2}e_{3}\right) .  \label{bhh4}
\end{equation}%
Here we can assume that $\alpha _{2}^{\prime }\neq 0$ without loss of
generality (when $\alpha _{2}^{\prime }=0$ then $\gamma $ becomes a geodesic
again). Then from (\ref{bhh4}) it follows that%
\begin{eqnarray}
k_{1}^{\,2}\varepsilon _{2} &=&g(\nabla _{T}T,\nabla _{T}T)  \nonumber \\
&=&\left( \alpha _{2}^{\prime }\right) ^{2}\left( \cosh ^{2}\alpha
_{2}-\sinh ^{2}\alpha _{2}\right)  \nonumber \\
&=&\left( \alpha _{2}^{\prime }\right) ^{2}.  \label{bhh4a}
\end{eqnarray}%
If $N$ is timelike then $k_{1}=0$ which is a contradiction again. This
completes the proof.

\begin{prop}
Let $\gamma :I\rightarrow {\cal HH}_{3}$ be a non-geodesic non-null curve
parametrized by arclenght with $B_{3}=0.$ Then $k_{2}^{2}=1$ and $\gamma $
cannot be biharmonic.
\end{prop}

\noindent {\em Proof.} Assume that $\gamma :I\rightarrow {\cal HH}_{3}$ is a
non-geodesic non-null curve parametrized by arclenght and $\gamma ^{\,\prime
}(s)=T(s)$. If $\gamma $ is a spacelike curve then from Proposition \ref%
{prop-bhh1} and (\ref{bhh2a}), $N$ must be timelike and $k_{1}=\,\pm \alpha
_{1}^{\prime }\neq 0.$ Using (\ref{bhh1}), (\ref{bhh2}), the first Frenet
equation and the definition of cross product in ${\cal HH}_{3}$ it follows
that%
\begin{eqnarray*}
N &=&\mp \left( \sinh \alpha _{1}\cosh \beta _{1}e_{1}+\sinh \alpha
_{1}\sinh \beta _{1}e_{2}+\cosh \alpha _{1}e_{3}\right) , \\
B &=&T\times N=\pm \left( \sinh \beta _{1}e_{1}+\cosh \beta _{1}e_{2}\right)
.
\end{eqnarray*}%
From (\ref{hh4}) we also have
\begin{eqnarray*}
\nabla _{T}N &=&\mp \lbrack \left( \alpha _{1}^{\prime }\cosh \alpha
_{1}\cosh \beta _{1}-\sinh \beta _{1}\right) e_{1} \\
&&+\left( \alpha _{1}^{\prime }\cosh \alpha _{1}\sinh \beta _{1}-\cosh \beta
_{1}\right) e_{2} \\
&&+\alpha _{1}^{\prime }\sinh \alpha _{1}e_{3}].
\end{eqnarray*}%
which implies that%
\[
k_{2}=g(\nabla _{T}N,B)=-1.
\]

Similarly, if $\gamma $ is a timelike curve then from Proposition \ref%
{prop-bhh1} and (\ref{bhh4a}), we have $N$ is spacelike and $k_{1}=\,\pm
\alpha _{2}^{\prime }\neq 0.$ Using (\ref{bhh4}) and the first Frenet
equation one obtains
\begin{eqnarray*}
N &=&\pm \left( \cosh \alpha _{2}\cosh \beta _{2}e_{1}+\cosh \alpha
_{2}\sinh \beta _{2}e_{2}+\sinh \alpha _{2}e_{3}\right) , \\
B &=&T\times N=\pm \left( \sinh \beta _{2}e_{1}+\cosh \beta _{2}e_{2}\right)
\end{eqnarray*}%
After a straightforward computation we get
\begin{eqnarray*}
\nabla _{T}N &=&\pm \lbrack \left( \alpha _{2}^{\prime }\sinh \alpha
_{2}\cosh \beta _{2}+\sinh \beta _{2}\right) e_{1} \\
&&+\left( \alpha _{2}^{\prime }\sinh \alpha _{2}\sinh \beta _{2}+\cosh \beta
_{2}\right) e_{2} \\
&&+\alpha _{2}^{\prime }\cosh \alpha _{2}e_{3}].
\end{eqnarray*}%
which gives%
\[
k_{2}=g(\nabla _{T}N,B)=-1.
\]%
The proof is completed.

\bigskip

Thus we have

\begin{cor}
Let $\gamma :I\rightarrow {\cal HH}_{3}$ be a non-geodesic non-null
biharmonic helix parametrized by arclenght. Then
\begin{equation}
\left\{
\begin{array}{c}
B_{3}=\text{constant}\neq 0, \\
k_{1}^{2}\varepsilon _{1}\varepsilon _{3}+k_{2}^{2}=1+4\varepsilon
_{3}B_{3}^{2}, \\
N_{3}=0.%
\end{array}%
\right.  \label{bhh5}
\end{equation}

\begin{lem}
Let $\gamma :I\rightarrow {\cal HH}_{3}$ be a non-geodesic non-null curve
parametrized by arclenght. If $N_{3}=0$ then
\begin{equation}
T(s)=\cosh \alpha _{0}\cosh \beta (s)e_{1}+\cosh \alpha _{0}\sinh \beta
(s)e_{2}+\sinh \alpha _{0}e_{3}  \label{bhh6}
\end{equation}%
or%
\begin{equation}
T(s)=\sinh \nu _{0}\cosh \rho (s)e_{1}+\sinh \nu _{0}\sinh \rho
(s)e_{2}+\cosh \nu _{0}e_{3},  \label{bhh7}
\end{equation}%
where $\alpha _{0}$, $\nu _{0}\in R$.
\end{lem}
\end{cor}

\noindent {\em Proof. }Let $T$ be the tangent vector field of $\gamma
:I\rightarrow {\cal HH}_{3}$ given by $T=T_{1}e_{1}+T_{2}e_{2}+T_{3}e_{3}$
and $g(T,T)=\varepsilon _{1}$. By using (\ref{hh4}) we have
\begin{eqnarray*}
\nabla _{T}T &=&(T_{1}^{\,\prime }-2T_{2}T_{3})e_{1}+(T_{2}^{\,\prime
}-2T_{1}T_{3})e_{2}+T_{3}^{\,\prime }e_{3} \\
&=&k_{1}\varepsilon _{2}N,
\end{eqnarray*}%
which implies that $N_{3}=0$ if and only if $T_{3}=$ constant. Then we
complete the proof.

\begin{th}
\label{theorem-parametric-spacelike}The parametric equations of all
non-geodesic spacelike biharmonic curves $\gamma $ of ${\cal HH}_{3}$ are%
\begin{eqnarray}
x(s) &=&\frac{1}{a}\cosh \alpha _{0}\sinh \left( as+b\right) +c_{1},
\nonumber \\
y(s) &=&\frac{1}{a}\cosh \alpha _{0}\cosh \left( as+b\right) +c_{2},
\label{bhh8} \\
z(s) &=&2\left( \sinh \alpha _{0}-\frac{1}{a}\left( \cosh \alpha _{0}\right)
^{2}\right) s  \nonumber \\
&&+\frac{2c_{1}}{a}\cosh \alpha _{0}\cosh \left( as+b\right) -\frac{2c_{2}}{a%
}\cosh \alpha _{0}\sinh \left( as+b\right) +c_{3},  \nonumber
\end{eqnarray}%
where $a=\sinh \alpha _{0}\pm \sqrt{5\left( \sinh \alpha _{0}\right) ^{2}+1}%
, $ $b,$ $c_{i}$ $\in R,$ $\left( 1\leq i\leq 3\right) $.
\end{th}

\noindent {\em Proof. }Assume that $\gamma :I\rightarrow {\cal HH}_{3}$ be a
spacelike non-geodesic curve. Then its tangent vector field is given by (\ref%
{bhh6}). From Gram-Schmidt procedure we have%
\[
N(s)=\sinh \beta (s)e_{1}+\cosh \beta (s)e_{2}.
\]%
By taking covariant derivative of the vector field $T$ we get%
\begin{eqnarray*}
\nabla _{T}T &=&\cosh \alpha _{0}\left( \beta ^{\,\prime }-2\sinh \alpha
_{0}\right) \left( \sinh \beta \,e_{1}+\cosh \beta \,e_{2}\right) \\
&=&k_{1}\varepsilon _{2}N,
\end{eqnarray*}%
where
\begin{equation}
k_{1}=\left\vert \cosh \alpha _{0}\left( \beta ^{\,\prime }-2\sinh \alpha
_{0}\right) \right\vert .  \label{bhh10}
\end{equation}%
Taking into account the cross product in ${\cal HH}_{3}$ one obtains%
\begin{eqnarray}
B(s) &=&T(s)\times N(s)  \nonumber \\
&=&\sinh \alpha _{0}\cosh \beta (s)e_{1}+\sinh \alpha _{0}\sinh \beta
(s)e_{2}+\cosh \alpha _{0}e_{3}.  \label{bhh11}
\end{eqnarray}%
Moreover,%
\[
\nabla _{T}N=\cosh \beta \left( \beta ^{\,\prime }-\sinh \alpha _{0}\right)
e_{1}+\sinh \beta \left( \beta ^{\,\prime }-\sinh \alpha _{0}\right)
e_{2}+\cosh \alpha _{0}e_{3}.
\]%
From the second Frenet equation, it follows that%
\begin{equation}
k_{2}=\sinh \alpha _{0}\left( \beta ^{\,\prime }-2\sinh \alpha _{0}\right)
-1.  \label{bhh12}
\end{equation}%
Then $\gamma $ is a spacelike non-geodesic biharmonic curve if and only if
\begin{equation}
\left\{
\begin{array}{c}
\beta ^{\,\prime }=\text{constant}\neq 2\sinh \alpha _{0}, \\
-k_{1}^{2}+k_{2}^{2}=1-4B_{3}^{2}.%
\end{array}%
\right.  \label{bhh13}
\end{equation}%
By substituting (\ref{bhh10}), (\ref{bhh12}) and $B_{3}=\cosh \alpha _{0}$
in the second equation of (\ref{bhh13}) we get%
\[
\left( \beta ^{\,\prime }\right) ^{2}-2\beta ^{\,\prime }\left( \sinh \alpha
_{0}\right) -4-4\left( \sinh \alpha _{0}\right) ^{2}=0
\]%
which gives%
\[
\beta ^{\,\prime }=\sinh \alpha _{0}\pm \sqrt{5\left( \sinh \alpha
_{0}\right) ^{2}+1}=a,
\]%
that is,%
\[
\beta (s)=as+b,\quad b\in R.
\]%
To find a differential equation system for the non-geodesic spacelike
biharmonic curve $\gamma (s)=\left( x(s),y(s),z(s)\right) ,$ by using (\ref%
{hh3}) we first note that
\begin{equation}
\frac{\partial }{\partial x}=e_{1}+ye_{3},\quad \frac{\partial }{\partial y}%
=e_{2}-xe_{3},\quad \frac{\partial }{\partial z}=\frac{1}{2}e_{3}.
\label{bhh13a}
\end{equation}%
Therefore since $T=\frac{d\gamma }{ds}$, we have the following differential
equations system%
\[
\begin{array}{c}
\frac{dx}{ds}=\cosh \alpha _{0}\cosh \left( as+b\right) , \\
\frac{dy}{ds}=\cosh \alpha _{0}\sinh \left( as+b\right) , \\
\frac{dz}{ds}=2\sinh \alpha _{0}+2\cosh \alpha _{0}\left( \sinh \left(
as+b\right) x(s)-\cosh \left( as+b\right) y(s)\right) .%
\end{array}%
\]%
Integrating the system gives (\ref{bhh8}). The proof is completed.

\begin{th}
\label{theorem-parametric-timelike}The parametric equations of all
non-geodesic timelike biharmonic curves $\gamma $ in ${\cal HH}_{3}$ are%
\begin{eqnarray}
\widetilde{x}(s) &=&\frac{1}{\widetilde{a}}\sinh \nu _{0}\sinh \left(
\widetilde{a}s+\widetilde{b}\right) +d_{1},  \nonumber \\
\widetilde{y}(s) &=&\frac{1}{\widetilde{a}}\sinh \nu _{0}\cosh \left(
\widetilde{a}s+\widetilde{b}\right) +d_{2},  \label{bhh14} \\
\widetilde{z}(s) &=&2\left( \cosh \nu _{0}-\frac{1}{\widetilde{a}}\left(
\sinh \nu _{0}\right) ^{2}\right) s  \nonumber \\
&&+\frac{2d_{1}}{\widetilde{a}}\sinh \nu _{0}\cosh \left( \widetilde{a}s+%
\widetilde{b}\right) -\frac{2d_{2}}{\widetilde{a}}\sinh \nu _{0}\sinh \left(
\widetilde{a}s+\widetilde{b}\right) +d_{3},  \nonumber
\end{eqnarray}%
where $\widetilde{a}=\cosh \nu _{0}\pm \sqrt{5\left( \cosh \nu _{0}\right)
^{2}-1},$ $\widetilde{b},$ $d_{i}$ $\in R,$ $\left( 1\leq i\leq 3\right) $.
\end{th}

\noindent {\em Proof. }The tangent vector field of a non-geodesic timelike
biharmonic curve $\gamma :I\rightarrow {\cal HH}_{3}$ can be given by (\ref%
{bhh7}). From Gram-Schmidt procedure we have%
\[
N(s)=\sinh \rho (s)e_{1}+\cosh \rho (s)e_{2},
\]%
which implies that $N$ is a timelike vector field. If we take the covariant
derivative of the tangent vector field $T$ it is easy to see that
\begin{eqnarray*}
\nabla _{T}T &=&\sinh \nu _{0}\left( \rho ^{\,\prime }-2\cosh \nu
_{0}\right) \left( \sinh \rho \,e_{1}+\cosh \rho \,e_{2}\right) \\
&=&k_{1}\varepsilon _{2}N
\end{eqnarray*}%
and
\begin{equation}
k_{1}=\left\vert \sinh \nu _{0}\left( \rho ^{\,\prime }-2\cosh \nu
_{0}\right) \right\vert .  \label{bhh15}
\end{equation}%
Also we have%
\begin{eqnarray}
B(s) &=&T(s)\times N(s)  \nonumber \\
&=&\cosh \nu _{0}\cosh \rho (s)e_{1}+\cosh \nu _{0}\sinh \rho (s)e_{2}+\sinh
\nu _{0}e_{3}.  \label{bhh16}
\end{eqnarray}%
In this case it is obvious that $B$ is a spacelike vector field. From (\ref%
{hh4}) we get%
\[
\nabla _{T}N=\cosh \rho \left( \rho ^{\,\prime }-\cosh \nu _{0}\right)
e_{1}+\sinh \rho \left( \rho ^{\,\prime }-\cosh \alpha _{0}\right)
e_{2}+\sinh \nu _{0}e_{3}.
\]%
It follows that%
\begin{equation}
k_{2}=\cosh \nu _{0}\left( \beta ^{\,\prime }-2\cosh \nu _{0}\right) +1.
\label{bhh17}
\end{equation}%
Then $\gamma $ is biharmonic if and only if
\begin{equation}
\left\{
\begin{array}{c}
\rho ^{\,\prime }=\text{constant}\neq 2\cosh \nu _{0}, \\
-k_{1}^{2}+k_{2}^{2}=1+4B_{3}^{2}.%
\end{array}%
\right.  \label{bhh18}
\end{equation}%
Using (\ref{bhh15}), (\ref{bhh17}) and $B_{3}=\sinh \nu _{0}$ in the second
equation of (\ref{bhh18}) we get%
\[
\left( \rho ^{\,\prime }\right) ^{2}-2\rho ^{\,\prime }\left( \cosh \nu
_{0}\right) +4-4\left( \cosh \nu _{0}\right) ^{2}=0
\]%
which gives%
\[
\rho ^{\,\prime }=\cosh \nu _{0}\pm \sqrt{5\left( \cosh \nu _{0}\right)
^{2}-1}=\widetilde{a},
\]%
that is,%
\[
\rho (s)=\widetilde{a}s+\widetilde{b},\quad \widetilde{b}\in R.
\]%
Since $T=\frac{d\gamma }{ds}$, from (\ref{bhh13a}), the differential
equations system for the non-geodesic timelike biharmonic curve $\gamma
(s)=\left( \widetilde{x}(s),\widetilde{y}(s),\widetilde{z}(s)\right) $ is
the following

\[
\begin{array}{c}
\frac{d\widetilde{x}}{ds}=\sinh \nu _{0}\cosh \left( \widetilde{a}s+%
\widetilde{b}\right) , \\
\frac{d\widetilde{y}}{ds}=\sinh \nu _{0}\cosh \left( \widetilde{a}s+%
\widetilde{b}\right) , \\
\frac{d\widetilde{z}}{ds}=2\cosh \nu _{0}+2\sinh \nu _{0}\left( \sinh \left(
\widetilde{a}s+\widetilde{b}\right) \widetilde{x}(s)-\cosh \left( \widetilde{%
a}s+\widetilde{b}\right) \widetilde{y}(s)\right) .%
\end{array}%
\]%
If we integrate the above system gives (\ref{bhh14}).

From Theorem \ref{theorem-parametric-spacelike} and Theorem \ref%
{theorem-parametric-timelike} we also have

\begin{cor}
\label{corollaryN3=0}Let $\gamma :I\rightarrow {\cal HH}_{3}$ be a
non-geodesic non-null curve parametrized by arclenght with $N_{3}=0.$ Then
we have $\varepsilon _{1}=-\varepsilon _{3}$ and $N$ is a timelike vector
field, where $\varepsilon _{1}=g(T,T)$ and $\varepsilon _{3}=g(B,B)$.
\end{cor}

\section{Horizontal Biharmonic curves in $3$-dimensional\newline
hyperbolic Heisenberg group}

Let $(x,y)\longrightarrow H_{(x,y)}$ be a non-integrable two dimensional
distribution in $R^{3}=R_{(x,y)}^{2}\times R$ defined by $H=\ker w$. The
distribution $H$ is said to be the horizontal distribution. A curve $%
s\longrightarrow \gamma (s)$, $\gamma (s)=\left( x(s),y(s),z(s)\right) $ is
called horizontal curve if $\gamma ^{\prime }(s)\in H_{\gamma (s)},$ for all
$s.$ By using (\ref{bhh13a}), for a non-null curve $\gamma $ in $3$%
-dimensional hyperbolic Heisenberg group we can write
\begin{equation}
\gamma ^{\prime }(s)=x^{\prime }(s)\frac{\partial }{\partial x}+y^{\prime
}(s)\frac{\partial }{\partial y}+z^{\prime }(s)\frac{\partial }{\partial z}%
=x^{\prime }(s)e_{1}+y^{\prime }(s)e_{2}+w(\gamma ^{\prime }(s))\frac{%
\partial }{\partial z}.  \label{hbc-1}
\end{equation}%
Then $\gamma $ is a horizontal curve if%
\begin{equation}
\gamma ^{\prime }(s)=x^{\prime }(s)e_{1}+y^{\prime }(s)e_{2},  \label{hbc-2}
\end{equation}
\begin{equation}
w(\gamma ^{\prime }(s))=z^{\prime }(s)+2x^{\prime }(s)y(s)-2x(s)y^{\prime
}(s).  \label{hbc-3}
\end{equation}

\begin{th}
The parametric equations of all non-geodesic spacelike horizontal biharmonic
curves $\gamma$ in ${\cal HH}_{3}$ are%
\begin{eqnarray}
x(s) &=&\pm \sinh \left( \pm s+b\right) +c_{1},  \nonumber \\
y(s) &=&\pm \cosh \left( \pm s+b\right) +c_{2},  \label{hbc-4} \\
z(s) &=&\mp 2s\pm 2c_{1}\cosh \left( \pm s+b\right) \mp 2c_{2}\sinh \left(
\pm s+b\right) +c_{3},  \nonumber
\end{eqnarray}%
where $b,$ $c_{i}$ $\in R,$ $\left( 1\leq i\leq 3\right) .$
\end{th}

\noindent {\em Proof. }Let $\gamma :I\rightarrow {\cal HH}_{3}$ be a
non-geodesic spacelike horizontal biharmonic curve. Since the tangent vector
field of $\gamma $ can be written as $T=T_{1}e_{1}+T_{2}e_{2}+T_{3}e_{3}$
then from (\ref{bhh6}) and (\ref{hbc-2}) we have
\begin{equation}
T_{3}=\sinh \alpha _{0}=0.  \label{hbc-5}
\end{equation}%
By using the last equation in (\ref{bhh8}) we complete the proof.

\begin{th}
There does not exist a non-geodesic timelike horizontal biharmonic curve in $%
{\cal H}_{3}.$
\end{th}

\noindent {\em Proof.} Assume that $\gamma :I\rightarrow {\cal HH}_{3}$ is a
non-geodesic timelike horizontal biharmonic curve. Then we have $N_{3}=0$
and $T_{3}=0.$ Since $\gamma $ is a timelike curve then Corollary \ref%
{corollaryN3=0} implies that $N$ is a timelike and $B$ is a spacelike vector
field. Using (\ref{hh4}) we have
\begin{equation}
g(\nabla _{T}T,e_{3})=T_{3}^{\,\prime },\quad g(\nabla
_{T}N,e_{3})=N_{3}^{\,\prime }-T_{2}N_{1}+T_{1}N_{2},\quad g(\nabla
_{T}B,e_{3})=B_{3}^{\,\prime }-T_{2}B_{1}+T_{1}B_{2}.  \label{hbc-6}
\end{equation}%
On the other hand from the Frenet formulas one can easily see that
\begin{subequations}
\begin{equation}
g(\nabla _{T}T,e_{3})=-k_{1}N_{3},\quad g(\nabla
_{T}N,e_{3})=k_{1}T_{3}+k_{2}B_{3},\quad g(\nabla _{T}B,e_{3})=-k_{2}N_{3}.
\label{hbc-7}
\end{equation}%
It follows from the definition of the cross product in ${\cal HH}_{3}$, (\ref%
{hbc-6}) and (\ref{hbc-7}) that
\end{subequations}
\[
k_{2}=1.
\]%
Substituting the last equation in (\ref{bhc-14}) we get
\[
-k_{1}^{2}=4B_{3}^{2},
\]%
which is a contradiction. The proof is completed.

{\it Authors' address}:\newline
Selcen Y\"{U}KSEL PERKTA\c{S} and Erol KILI\c{C},\newline
Department of Mathematics,\newline
Inonu University, 44280, Malatya/TURKEY\newline
E-mail: selcenyuksel@@inonu.edu.tr, ekilic@@inonu.edu.tr

\bigskip

\end{document}